\documentclass[12pt]{article}
\usepackage{cite}
\usepackage{amssymb}
\usepackage{amsmath}
\numberwithin{equation}{section}

\newcommand{\eh}{\hfill}\newlength{\sperr}

\newenvironment{proof}{{\settowidth{\sperr}{\bf\rm
Proof}%
\par\addvspace{0.3cm}\noindent\parbox[t]{1.3\sperr}
{\bf\rm P\eh r\eh o\eh o\eh f.\eh }%
}}{\nopagebreak\mbox{}\hfill
$\Box$\par\addvspace{0.3cm}}

\def\Vol{\mathrm{Vol}}

\def\dv{{\mathrm{div}}}

\def\p{\partial}
\def\Om{\Omega}
\def\ve{\varepsilon}
\def\nn{\nonumber}

\def\BR{{\mathbb R}}

\def\wt{\widetilde}
\newcommand{\E}{\mathrm{e}}
\newtheorem{Pa}{Paper}[section]
\newtheorem{Tm}[Pa]{{\bf Theorem}}
\newtheorem{La}[Pa]{{\bf Lemma}}
\newtheorem{Cy}[Pa]{{\bf Corollary}}
\newtheorem{Rk}[Pa]{{\bf Remark}}

\newtheorem{Dn}[Pa]{{\bf Definition}}
\newtheorem{Nn}[Pa]{{\bf Notation}}
\newtheorem{Pn}[Pa]{{\bf Proposition}}
\title{Nonlinear Fokker-Planck  equation: stability, distance and the corresponding extremal problem
in the spatially inhomogeneous case}
\author{Alexander Sakhnovich, Lev Sakhnovich}
\date{}
\begin{document}
\maketitle
\thanks{A.L. Sakhnovich,   Fakult\"at f\"ur Mathematik,
Universit\"at Wien,
\\
Nordbergstrasse 15, A-1090 Wien, Austria, \\
E-mail: {\tt al$_-$sakhnov@yahoo.com }}\\

\thanks{L.A. Sakhnovich, 99 Cove ave., Milford, CT, 06461, USA, \\
 E-mail: lsakhnovich@gmail.com}\\


\begin{abstract}
We start with a  global Maxwellian
 $M_{k}$, which is a stationary solution, with the constant total density ($\rho(t)\equiv \wt \rho$), of the Fokker-Planck equation.
The notion of distance between the
 function $M_{k}$ and an arbitrary solution $f$
(with the same  total density $\wt \rho$ at the fixed moment $t$) of the Fokker-Planck equation is introduced.
In this way, we essentially generalize the important Kullback-Leibler distance, which was studied
before. Using this generalization, we show local stability of the global Maxwellians
 in the spatially inhomogeneous case. We compare also the energy and entropy in the classical and quantum cases.
\end{abstract}

\section{Introduction}
We  consider the  Fokker-Planck
 equation
 \begin{equation}\label{1}
 \frac{\partial{f}}{\partial{t}}=\Delta_{v}{f}-v\cdot \bigtriangledown_{x}{f}+
 \dv_{v}\big(vf(1+kf)\big),\end{equation}
 where  $t{\in}\BR$ stands for time,   $x=(x_{1},x_{2},...,x_{n}){\in}\Omega$ stands for
 space coordinates,  $v=(v_{1},v_{2},...,v_{n})\, {\in} \, \BR^{n}$ is velocity
 and $\, \BR \,$ denotes the real axis. This non-linear Fokker-Planck equation serves  as a kinetic model for bosons $(k>0)$ and fermions $(k<0)$.

The important notion of Kullback-Leibler distance \cite{KL} is essentially generalized in our paper and new
conditional extremal problems, which appear in this way, are solved. The solutions $f(t,x,v)$
of the Fokker-Planck equation are studied in the bounded domain
$\Omega$ of the $x$-space. Such an approach essentially changes the usual situation, that is,
the total energy depends on $t$ and
the notion of distance (between
a stationary solution and an  arbitrary solution of the Fokker-Planck equation)
 includes the $x$-space.  Thus,  the notion of distance remains well-defined also in the spatially inhomogeneous case. 
 Recall that the Kullback-Leibler distance, which has many applications (see, e.g., \cite{Hab, SH, Vi2} and references therein), 
 is  defined only in the spatially homogeneous case.
 
 In our previous paper \cite{ASLS} we studied a model case of the one dimensional $x$-space.
 Here the case $\dim \Omega \geq 1$ is dealt with.
Furthermore, using our generalization of the Kullback-Leibler distance,  we show 
local stability of global Maxwellians
 in the spatially inhomogeneous case. 
 
 The comparison of the energy and entropy in the classical and quantum cases is an important subject
(see  \cite{Caret, LAS1, LAS2, LAS3, LASPLet, LAS4, We} and references therein).
 Here, we compare these energy and entropy for the situation described by the Fokker-Planck
 equation. It is especially interesting for the applications that the fermion and  boson cases are essentially different.

Our definition  of the quantum entropy $S_k$ ($k\not=0$) is slightly different from the previous definitions (see \cite{Dol, Lu2}).
We  show that the natural requirement
\begin{equation}\label{q1}
S_{k}{\to}S_{c},\quad k{\to}0 \quad (S_c=S_0 \,\,\mathrm{is} \,\, \mathrm{the} \,\,
\mathrm{classical} \,\, \mathrm{entropy})
\end{equation}
is not fulfilled in the case of the old definition, however \eqref{q1}  holds for   our modified definition (see Section \ref{Se3}).
Some necessary definitions are given in Subsection \ref{Prel}.
An important functional, which attains its maximum at the
function $M_{k}$ is introduced there.
The distance between solutions and the corresponding extremal problem
are studied in Sections \ref{Se3} and \ref{Lya}.
Our results on the Fokker-Planck equation are mostly related to the corresponding results
 (from \cite{LASPLet, LAS4})
on the Boltzmann equation but the theorems from Section \ref{Lya}
have no analogs in the Boltzmann case.

We use the standard notation $|v|=\sqrt{v_1^2+\ldots + v_n^2}$ and $C_{0}^{1}$  denotes the class
of differentiable functions $f(x,v)$, which tend to zero sufficiently rapidly
when $v$ tends to infinity.
\section{The extremal problem} \label{Se3}
 \subsection{Preliminaries}\label{Prel}
Here, we  present some  well-known  notions and results
 connected with   equation \eqref{1}.
It is required that the distribution function $f(t,x,v)$ satisfies the inequalities
\begin{equation}\label{3} 
f(t,x,v) \geq 0, \quad 1+kf(t,x,v) \geq 0 \quad (k \in \BR),
\end{equation}
and we set $f\log f=0$ for the case that $f=0$ and $(1+kf)\log(1+kf)=0$ for $1+kf=0$.
Then, the mapping
\begin{align}\label{5'}&
\Phi(f):=f\log{f}\quad {\mathrm{for}} \quad k=0, \\
\label{5}&
 \Phi(f):=f\log{f}-\frac{1}{k}(1+kf)\log(1+kf)+f \quad {\mathrm{for}} \quad k\not=0, 
 \end{align}
 is well-defined. Now, the entropy is given by the equality
 \begin{equation}\label{4}
 S(t,f)=S(t)=-\int_{\Omega}\int_{\BR^n}\Phi(f)d{v}dx.
\end{equation}
The notions of density $\rho(t,x)$,  total density $\rho(t)$, mean velocity $u(t,x)$,
energy $E(t,x)$, and total  energy $E(t)$ are introduced via formulas:
 \begin{align}\label{7}&
 \rho(t,x)=\int_{\BR^n} f(t,x,v)dv , \quad \rho(t)=\int_{\Omega}\rho(t,x){dx},\\
 \label{8}&
 u(t,x)=\big(1/\rho(x,t)\big)\int_{\BR^n} v f (t,x,v){d}v,
 \\ \label{9}&
 E(t,x)=\int_{\BR^n} \frac{|v|^{2}}{2}f(t,x,v)dv , \quad
E(t)=\int_{\Omega}\int_{\BR^n} \frac{|v|^{2}}{2}f(t,x,v)dv dx.
\end{align}
We assume that the domain $\Omega$ is bounded, and so its volume
is also bounded:
\begin{equation}\label{11}
\Vol(\Omega)=V_{\Omega}<\infty.
\end{equation}
 \subsection{The free energy functional and extremal problem}\label{Free}
We  introduce the "free energy" functional
\begin{equation} \label{12}
F(f)=F\big(f(t)\big)=S(t) - E(t),
\end{equation}
 where  $S(t)$ and $E(t)$ are defined by formulas \eqref{4} and \eqref{9}, respectively.
Next,   we use the calculus
of variations (see \cite{Ha}) and find the function $f_{max}$ which maximizes
the functional \eqref{12}, where the parameters  $\, t\,$ and $\rho(t)=\wt \rho>0$ are fixed.
 The corresponding Euler's equation takes the form
 \begin{equation} \label{14}
 -\frac{|v|^{2}}{2}- \log{f}+\log(1+kf)+
{\mu}=0.
\end{equation}
From the last relation we obtain
 \begin{equation} \label{15}
 f/(1+kf)=C\exp\left\{-{|v|^{2}}/({2})\right\}, \quad C:=\E^{-\mu}.
 \end{equation}
 Formula \eqref{15} implies  that
 \begin{equation} \label{16}  f=M_{k}=\frac{C\exp\left\{-{|v|^{2}}/2\right\}}
{1-kC\exp\left\{-{|v|^{2}}/2\right\}},
\end{equation}
that is, $f$ coincides with the global Maxwellian $M_k$. 

In view of the requirement $\rho(t)=\wt \rho$, the constant $C$ in the equality \eqref{16} is derived from the relation
 \begin{equation} \label{16'}
 V_{\Omega}\int_{\BR^n}\frac{C\exp\left\{-{|v|^{2}}/2\right\}}
{1-kC\exp\left\{-{|v|^{2}}/2\right\}}dv=\wt \rho.
\end{equation}
The function  $f$ given by \eqref{15} (or, equivalently, by \eqref{16}) is nonsingular and satisfies conditions \eqref{3} and $\rho(t)=\wt \rho>0$
if and only if
\begin{equation}\label{op1}
C>0, \quad 1-kC>0.
\end{equation}
In Subsection \ref{comp} we prove that there is a unique value $C$ satisfying  relations
\eqref{16'} and \eqref{op1}.  Let us show that $F$ indeed attains  its maximum on the global Maxwellian $M_k$ 
corresponding to such $C$.
According to \eqref{4}, \eqref{9} and \eqref{12}  the "free energy" $F$ admits the representation
\begin{equation}\label{op2}
F=\int_{\Om}\int_{\BR^n}\Psi(f)dvdx, \quad \Psi:=-\frac{|v|^2}{2}f-\Phi.
\end{equation}
Taking into account \eqref{5'}, \eqref{5}, \eqref{op1} and \eqref{op2}, we have the inequality
 \begin{equation} \label{17}  \frac{\delta^{2}}{\delta{f}^{2}}\Psi=-\frac{1}{(1+kf)f}<0 ,
 \end{equation}
and the next proposition follows.
\begin{Pn}\label{Pn1} Under condition \eqref{op1}, the functional  $F$ given by \eqref{12}
attains its maximum on the function $M_k$ of the form \eqref{16}
$($where $C$ is defined in \eqref{16'}$)$, that is, 
 \begin{equation} \label{18}
 G(f)=F(M_{k})-F(f)>0\quad (f{\ne}M_{k}).
\end{equation}
\end{Pn}
\begin{Rk} In this subsection we introduced the important functional $F$, which attains its maximum
on the global Maxwellians and the generalization $G$ $($see \eqref{18}$)$ of the Kullback-Leibler distance.
As opposed to the Kullback-Leibler distance, which is defined in the $x$-homogeneous case, the distance $G$ is well-defined for the functions $f$,
which depend on $x$.

Later in this section we will consider the energy $E$, entropy $S$ and free energy $F$ of the
global Maxwellians. The following sections are dedicated to the study of the general solutions
of the Fokker-Planck equation.
\end{Rk}
\subsection{Comparison of the classical and quantum characteristics}\label{comp}
 Let us calculate the integral on the left-hand side of \eqref{16'}. Using  spherical coordinates, we have
 \begin{equation} \label{3.9}
C \int_{\BR^n}\frac{\exp\left\{-{|v|^{2}}/2\right\}}
{1-kC\exp\left\{-{|v|^{2}}/2\right\}}dv=\omega_{n-1}C\int_{0}^{\infty}
\frac{r^{n-1}\exp\left\{-{r^{2}}/2\right\}}
{1-kC\exp\left\{-{r^{2}}/2\right\}}dr,
\end{equation}
where the surface area of the $(n-1)$-sphere of radius $1$ is
 \begin{equation} \label{3.10}
 \omega_{n-1}=\frac{2{\pi}^{n/2}}{\Gamma(n/2)},
 \end{equation}
and $\Gamma(z)$ is the gamma function. Euler's integral representation of the gamma function
easily yields
 \begin{equation} \label{3.11}
 \int_{0}^{\infty}\E^{-ar^{2}}r^{n-1}dr=\frac{1}{2}a^{-n/2}\Gamma(n/2) \quad (a>0).
 \end{equation}
 \begin{Nn}\label{NnC} According to \eqref{16} and \eqref{16'}, the value $C$  corresponding to $M_k$
 depends on $k$.  We denote this value by $C_k$.
 \end{Nn}
 Taking into account \eqref{16'} and \eqref{3.9}--\eqref{3.11},
we obtain
 \begin{equation} \label{3.12}
 (2{\pi})^{n/2}V_{\Omega}C_kL_{n/2}(kC_k)=\wt \rho,
\end{equation}
where
 \begin{equation} \label{3.13}
  L_{n/2}(z)=\frac{2^{1-(n/2)}}{\Gamma(n/2)}\int_{0}^{\infty}\frac{\E^{-\frac{r^{2}}{2}}}
{1-z\E^{-\frac{r^{2}}{2}}}r^{n-1}dr=\sum_{m=1}^{\infty}\frac{z^{m-1}}{m^{n/2}}.
\end{equation}
We note that the series representation $L_{n/2}(z)=\sum_{m=1}^{\infty}\frac{z^{m-1}}{m^{n/2}}$
does not hold for $|z|>1$, and we use only the first equality in \eqref{3.13} for the case that $z< -1$.
Using the first equality in  \eqref{3.13}, we derive the following statement.
\begin{Pn}\label{Proposition 3.1.}
The function $L_{n/2}(z)$  increases strictly monotonically in the interval $-\infty \leq z<1$ and
 \begin{align} \label{3.13'} &
L_{n/2}(0)=1; 
 \\ & \label{3.14}
 L_{1/2}(1)=L_{1}(1)=\infty;\quad L_{n/2}(1)<\infty,\,\, n>2.
 \end{align}
 \end{Pn}

Proposition \ref{Proposition 3.1.} implies the next two corollaries.

\begin{Cy}\label{CyB}
If $k>0$ $($boson case$)$ and either $n=1$ or $n=2$, then
equation \eqref{3.12} has one and only one solution $C_k$ such that $C_k>0$ and $kC_k<1$.
\end{Cy}

\begin{Cy}
\label{Cy 2.5}
If $k>0$ $($boson case$)$, $n>2$ and
 \begin{equation} \label{3.15}
 (2{\pi})^{n/2}V_{\Omega}L_{n/2}(1)>k\wt \rho,
 \end{equation}
then
equation \eqref{3.12} has one and only one solution $C_k$ such that $C_k>0$ and $kC_k<1$.
\end{Cy}

\begin{Rk}
The function $L_{n/2}(z)$ belongs to the class of the $L$-functions \cite{LaGa} and
is connected with the famous $($see, e.g., \cite{T}$)$ Riemann zeta-function 
 \begin{equation} \label{3.16}
 \zeta(z)=\sum_{p=1}^{\infty}\frac{1}{p^{z}};\quad \Re{z}>1
 \end{equation}
 by the relation
 \begin{equation} \label{3.17}
 L_{n/2}(1)=\zeta(n/2).
 \end{equation}
Hence, we have the equalities
 \begin{equation} \label{3.18}
 L_{3/2}(1)=\zeta(3/2)=2.612,\quad L_{2}(1)=\zeta(2)=1.645,
 \end{equation}
 \begin{equation} \label{3.19}
 L_{5/2}(1)=\zeta(5/2)=1.341,\quad L_{3}(1)=\zeta(3)=1.202.
\end{equation}
\end{Rk}

Let us study the fermion case $k<0$.  Taking into account Notation \ref{NnC},
we may later  (differently from the constants $C_k$) consider $C$ as a variable. 
In view of \eqref{3.13}, the next proposition is valid.

\begin{Pn}
\label{Proposition 3.2}
Assume that $k<0$. Then, $CL_{n/2}(kC)$ increases strictly monotonically with   respect to $C$ $\,(0 \leq C<\infty)$,
 and $CL_{n/2}(kC){\to}\infty$ for $C{\to}\infty$.
 \end{Pn}

\begin{Cy}\label{CyF}
If $k<0$ $($fermion case$)$, then
equation \eqref{3.12} has one and only one solution $C_k$ such that $C_k>0$.
\end{Cy}
The second inequality in \eqref{op1}  holds in the fermion case (i.e., in the case $\,\,C=C_k$ and $k<0$) automatically.

Finally,  we consider in this section
 the   energy of the global Maxwellian:
 \begin{align} \nn
 E_k=E(M_k)&=\frac{1}{2}\int_{\Omega} \int_{\BR^n}|v|^2\frac{C_k\exp\left\{-{|v|^{2}}/2\right\}}
{1-kC_k\exp\left\{-{|v|^{2}}/2\right\}}dvdx
\\ \label{3.20}
&=\frac{1}{2}\omega_{n-1}V_{\Om}C_k\int_{0}^{\infty}r^{n+1}
\frac{\exp\left\{-{r^{2}}/2\right\}}
{1-kC_k\exp\left\{-{r^{2}}/2\right\}}dr.
\end{align}
It is immediate from \eqref{3.12} that
 \begin{equation} \label{op3}
\wt \rho/ L_{n/2}(kC_k)=(2 \pi )^{n/2}V_{\Om}C_k,
\end{equation}
and so, using \eqref{3.10} and \eqref{3.13}, we rewrite \eqref{3.20} in the form
 \begin{equation} \label{3.21}
 E_{k}=\left(\frac{n{\wt \rho}}{2}\right)\frac{L_{(n/2)+1}(kC_k)}{L_{n/2}(kC_k)}.
\end{equation}
We note that the corresponding classical energy $E_c$ (i.e., the energy for the case $k=0$) is given by the formula
 \begin{equation} \label{3.22}
 E_{c}=\frac{n{\wt \rho}}{2}.
 \end{equation}
The points $kC=\pm 1$ are called the critical points in boson and fermion theories.
(Recall that the series representation of $L_{n/2}(z)$ in \eqref{3.13}
does not hold for $|z|>1$.)
 \begin{Pn}\label{Proposition 3.3.}
{Let the  condition}
 \begin{equation} \label{3.23}
-1\leq  kC_k{<}1
 \end{equation}
be fulfilled. Then, we have the inequalities
 \begin{equation} \label{3.24}
E_{q,B}<E_{c} <  E_{q,F}  \quad {\mathrm{for}} \quad n \geq 1, 
 \end{equation}
 where $E_{q}$ denotes the energy in the quantum case $($i.e., the case $k\not=0)$, $E_{q,B}$
 stands for the energy in the boson case $k>0$ and  $E_{q,F}$ stands for the energy in the fermion case $k<0$.
 \end{Pn}
\begin{proof}
Taking into account the second equality in \eqref{3.13}, we obtain:
 \begin{equation} \label{3.25}
  \frac{L_{(n/2)+1}(kC)}{L_{n/2}(kC)}<1\quad \mathrm{for}\quad k>0, \quad kC<1.
\end{equation}
Moreover, we will show that
 \begin{equation} \label{3.25'}
\frac{L_{(n/2)+1}(kC)}{L_{n/2}(kC)}>1\quad\mathrm{for} \quad k<0, \quad kC \geq -1.
\end{equation}
For this purpose, we compare sums of two consequent terms (with numbers $2p$ and $2p+1$) in the Taylor series representations
\eqref{3.13} of $L_{(n/2)+1}(z)$ and $L_{n/2}(z)$, and derive that
\begin{align}\nn &
 |z|^{2p-1}
 \left( 
 \left(-\frac{1}{(2p)^{l+1}}+\frac{|z|}{(2p+1)^{l+1}}\right)
 -\left(-\frac{1}{(2p)^{l}}+\frac{|z|}{(2p+1)^{l}}\right)
 \right)
 \\ & \label{op4}
 ={|z|^{2p-1}}\left(\frac{2p-1}{(2p)^{l+1}}-\frac{2p|z|}{(2p+1)^{l+1}}\right).
\end{align}
 Furthermore, it is immediate that
\begin{align}\label{op5} &
 {|z|^{2p-1}}\left(\frac{2p-1}{(2p)^{l+1}}-\frac{2p|z|}{(2p+1)^{l+1}}\right)=
 \frac{|z|^{2p-1}}{(2p)^l}\left(\frac{2p-1}{2p}-|z|\left(\frac{2p}{2p+1}\right)^{l+1}\right).
\end{align}
Finally, it is easy to see that, for $-1\leq z<0$, $l=n/2$, $n \geq 2$, we have
\begin{align}\label{op6} &
 \frac{2p-1}{2p}-|z|\left(\frac{2p}{2p+1}\right)^{l+1}\geq  \frac{2p-1}{2p}-\left(\frac{2p}{2p+1}\right)^{2}>0,
\end{align}
and relations \eqref{op4}--\eqref{op6} imply \eqref{3.25'}. It remains to prove \eqref{3.25'} for the case $n=1$.
We easily calculate directly that
\begin{align}\label{opp1} &
 \frac{2p-1}{2p}-\frac{8}{10}\left(\frac{2p}{2p+1}\right)^{3/2}>0,
\end{align}
which yields
\begin{align}\label{opp2} &
 \frac{2p-1}{2p}-|z|\left(\frac{2p}{2p+1}\right)^{3/2}>0 \quad {\mathrm{for}} \quad -0.8 \leq z <0.
\end{align}
Formulas \eqref{op4}, \eqref{op5} and \eqref{opp2} show that
\begin{align}\label{opp3} &
L_{3/2}(kC)>L_{1/2}(kC) \quad {\mathrm{for}} \quad -0.8 \leq kC <0.
\end{align}
We use connections between Lerch zeta functions and Riemann zeta functions and 
take into account the estimate \cite[sequence A078434]{oeis} in order
to calculate that up to the first two decimal places we have
\begin{align}\label{opp4} &
L_{3/2}(-1) =0.76 \, .
\end{align}
Taking into account that (for $-1\leq z<0$) the series in \eqref{3.13} is an alternating series
satisfying the Leibniz criterion, we obtain the inequalities
\begin{align}\label{opp5} &
0.65<L_{1/2}(-0.8)<0.6589 \, .
\end{align}
In view of Proposition \ref{Proposition 3.1.}, the functions $L_{1/2}$ and $L_{3/2}$ increase
monotonically on the interval $[-1, \,-0.8]$. Hence, relations \eqref{opp4} and \eqref{opp5} imply that
\begin{align}\label{opp6} &
L_{3/2}(kC)>L_{1/2}(kC) \quad {\mathrm{for}} \quad -1 \leq kC \leq -0.8 \, .
\end{align}
Inequalities \eqref{opp3} and \eqref{opp6} prove that \eqref{3.25'} holds also for $n=1$.
Thus, it is proved that \eqref{3.25'} is valid for all $n\geq 1$.
Inequalities \eqref{3.24} follow directly from \eqref{3.21}, \eqref{3.22}, \eqref{3.25} and \eqref{3.25'}.
\end{proof}
\begin{Rk}  The proof of formula \eqref{3.25'},  for the case that  $n=1$, 
shows  that Conjecture 6.1 from \cite{LASPLet} is valid.
\end{Rk}
Formulas \eqref{3.12} and \eqref{3.13'} yield
 \begin{equation} \label{op7}
 C_0=(2\pi)^{-n/2}\wt \rho /V_{\Om}.
 \end{equation}
\begin{La}\label{La2.10}
The following inequalities are valid:
 \begin{align} \label{3.30}&
 C_k>C_{0} \quad {\mathrm{for}} \quad k<0; \qquad  C_k<2C_{0} \quad {\mathrm{for}} \,\, -1<2kC_0<0;
 \\ \label{3.30'}& C_k<C_{0} \quad {\mathrm{for}} \quad 0<kC_0< 1.
 \end{align}
 \end{La}
 \begin{proof} 
 Let $k<0$. Then, according to Proposition \ref{Proposition 3.1.}, we have
\begin{equation} \label{op8}
 L_{n/2}(kC_0)<L_{n/2}(0)=1.
 \end{equation} 
 Using relations \eqref{3.12}, \eqref{op7} and \eqref{op8}, we obtain
 \begin{equation} \label{op9}
 C_kL_{n/2}(kC_k)=C_0>C_0L_{n/2}(kC_0).
 \end{equation}
 Hence, Proposition \ref{Proposition 3.2} implies that the first inequality in \eqref{3.30} holds.
 
 Next, let $-1<2kC_0<0$. Again using Proposition \ref{Proposition 3.1.}, we see that
  \begin{equation} \label{op10}
2 L_{n/2}(2kC_0)>2L_{n/2}(-1).
 \end{equation}
Taking into account \eqref{3.13} and \eqref{3.13'}, from \eqref{op10} we derive
 \begin{equation} \label{op11}
2 L_{n/2}(2kC_0)>L_{n/2}(0)=1.
 \end{equation}
 Thus, we obtain 
 \begin{equation} \label{op12}
 2 C_0 L_{n/2}(2kC_0)>C_0=(2\pi)^{-n/2}\wt \rho /V_{\Om}.
 \end{equation}
 In view of  \eqref{3.12} and \eqref{op12},  we have $2 C_0 L_{n/2}(2kC_0)>C_k L_{n/2}(kC_k)$.
 Then, as in the proof of the first inequality in \eqref{3.30}, we apply 
 Proposition \ref{Proposition 3.2} and see
  that  the second inequality in \eqref{3.30} holds.
 
 Finally, let $0<kC_0< 1$. Since $L_{n/2}(z)$ is increasing (see Proposition~\ref{Proposition 3.1.}),
 we have 
 \begin{equation} \label{op8'}
 L_{n/2}(kC_0)>L_{n/2}(0)=1,
 \end{equation} 
 and, moreover, $C L_{n/2}(kC)$ also increases strictly monotonically. According to \eqref{3.12}, \eqref{op7}
 and \eqref{op8'}, the relations
 \begin{equation} \label{op13}
 C_0L_{n/2}(kC_0)>C_0=C_kL_{n/2}(kC_k)
 \end{equation} 
 hold. Therefore, since $C L_{n/2}(kC)$ is monotonic, we see that  $C_k<C_0$.
 \end{proof}
 It follows from Lemma \ref{La2.10} that $C_k$ is bounded  in the  neighborhood of $k~=~0$.
 The behavior of $E_q$, $F_q=F(M_k)$ and the entropy $S_q=S(M_k)$ in the punctured neighborhood
 of  $k=0$ is given in the next proposition.
 \begin{Pn}\label{Pn3.4}
For $k\to 0$, we have  the following asymptotic  relations:
 \begin{align} \label{3.31}&
 E_{q}-E_{c}=-\frac{n{\wt \rho}kC_{0}}{{4}(2^{n/2})}+O(k^{2}) \quad (E_c=E(M_0)),
\\ & \label{3.32}
 S_{q}-S_{c}=-\frac{(n-2){\wt \rho}kC_{0}}{{4}(2^{n/2})}+O(k^{2}) \quad (S_c=S(M_0)),
\\ &  \label{3.33}
 F_{q}-F_{c}=\frac{{\wt \rho}kC_{0}}{{2}(2^{n/2})}+O(k^{2}) \quad (F_c=F(M_0)).
\end{align}
\end{Pn}
\begin{proof}
Step 1.
In order to calculate the entropy $S(M_{k})$ we recall definitions \eqref{4} and \eqref{16} of $S$ and $M_k$, respectively,
 and use the equalities
  \begin{equation} \label{n7}
M_{k}=g/(1-k g), \quad 1+k M_{k}=(1-k g)^{-1}, \quad g:=C_k\E^{-|v|^2/2},
\end{equation}
which simplify the expressions  $\Phi(M_k)$ for $\Phi$ given by \eqref{5'} and \eqref{5}:
\begin{align} \label{n5}&
\Phi(M_k)=M_{k}(1+\log g)+(1/k)\log(1-k g) \quad \mathrm{for} \quad k\not=0, \\
\label{n8}& 
\Phi(M_0)=g\log g \quad \mathrm{for}  \quad M_0=g. 
\end{align}
Recall definitions \eqref{7} and \eqref{9} of $\rho$ and $E$ and recall that for $f=M_k$ we have
$\rho(t) \equiv {\mathrm{const}}=\wt \rho$.
Substituting $\log g=\log C_k -(1/2)|v|^2$ into \eqref{n5}  and then substituting  \eqref{n5} 
 into
\eqref{4}, we obtain
 \begin{align} \label{n9}&
 S(M_{k})=E_{q}-(1+\log{C_k})\wt \rho -\frac{1}{k}V_{\Om}\int_{\BR^n}
 \log(1-k g)dv \quad  \mathrm{for} \quad k\not=0.
 \end{align}
Substituting $\log g=\log C_k -(1/2)|v|^2$ into \eqref{n8} and then substituting   \eqref{n8}
 into
\eqref{4}, we obtain
\begin{align}
  \label{3.26}& S(M_{0})=S_c= E_{c}-{\wt \rho}\log{C_{0}} ,
\end{align}
where $E_c=E(M_0)$ and $C_0$ is the value of $C_k$ for the case that $k=0$
(recall Notation \ref{NnC}). Taking into account the definition \eqref{9} of energy and 
using spherical coordinates and integration by parts,  we rewrite \eqref{n9}:
 \begin{equation} \label{3.27}
 S(M_{k})=E_{q}-(1+\log{C_k})\wt \rho 
+\frac{2}{n}{E_{q}} =\left(1+\frac{2}{n}\right){E_{q}}-(1+\log{C_k})\wt \rho \quad  \mathrm{for} \quad k\not=0.
\end{equation}

According to \eqref{3.22}, we have $\wt \rho=2E_c/n$. Therefore, for $k \not=0$ formulas \eqref{3.26} and \eqref{3.27} imply that
 \begin{equation} \label{3.28}
 S(M_k)-S(M_0)=S_{q}-S_{c}=\frac{n+2}{n}(E_{q}-E_{c})-{\wt \rho}\log({C_k}/{C_{0}}),
\end{equation}
Hence, taking into account  \eqref{12} and \eqref{3.28} we derive
 \begin{equation} \label{3.29}  F_{q}-F_{c}=\frac{2}{n}(E_{q}-E_{c})-{\wt \rho}\log({C_k}/{C_{0}}).
\end{equation}

Step 2. The equalities in  \eqref{op9} and \eqref{op13} yield
 \begin{equation} \label{op15}
\frac{kC_0}{L_{n/2}(kC_k)}=kC_k.
\end{equation} 
In view of formula \eqref{op7} and Lemma \ref{La2.10}, we see that the values
$|kC_k|$ and $\sup_{|z|\leq \ve}\left|\frac{d}{dz}\left(\frac{kC_0}{L_{n/2}(z)}\right)\right|$
are small 
for small values of $|k|$. 
Thus, we can apply the iteration method to the equation $z=\frac{kC_0}{L_{n/2}(z)}$
in order to derive
 \begin{equation} \label{op16}
C_k=C_0+O(k), \quad k \to 0.
\end{equation} 
Taking into account the series expansion in \eqref{3.13} and formulas \eqref{op15} and \eqref{op16},
we obtain
 \begin{equation} \label{op17}
C_k/C_0=1/{L_{n/2}(kC_k)}=1-kC_k/2^{n/2}+O(k^2).
\end{equation} 
Furthermore, from \eqref{op16} and \eqref{op17} we have
 \begin{equation} \label{op18}
\log\big(C_k/C_0\big)=-kC_0/2^{n/2}+O(k^2).
\end{equation}
Using formulas \eqref{3.21}, \eqref{3.22}, \eqref{op16} and the series expansion in \eqref{3.13},
we see that  \eqref{3.31} holds. According to \eqref{3.31} and \eqref{op18} we may rewrite \eqref{3.28}
in the form \eqref{3.32}. Finally, in view of \eqref{3.31} and \eqref{op18}, we rewrite \eqref{3.29}
in the form \eqref{3.33}.
\end{proof}

\begin{Cy}
The following inequalities hold for small values of $k:$
 \begin{equation} \label{3.34}
E_{q,B}<E_{c}<E_{q,F}; \quad S_{q,B}<S_{c}<S_{q,F} \quad {\mathrm{for}} \quad n>2; \quad F_{q,F}<F_{c}<F_{q,B}.
 \end{equation}
 \end{Cy}

\section{General-type solutions }\label{sec4}
\subsection{Dissipative and conservative solutions }
\paragraph{1.}
In this section we study  general solutions $f(t,x,v)$ (satisfying \eqref{3}) of the Focker-Planck equation \eqref{1}.
The total energy flux through the surface $\partial\Omega$ per unit time is given by the equalities
 \begin{align}\nn
A(f,\Omega): &= \int_{\Omega}\int_{\BR^n}(|v|^{2}/2)v{\cdot}{\bigtriangledown_{x}}f(t,x,v)d{v}dx
\\  \label{4.1}  &= \int_{\partial\Omega}\int_{\BR^n}(|v|^{2}/2)\big(v{\cdot}n(y)\big)f(t,y,v)d{v}d\sigma,
\end{align}
where $\partial\Omega$ is the boundary of the $\Omega$, and the integral $\int_{\partial\Omega}g d\sigma$ is the surface integral with $n(y)$ being the outward 
unit normal to that surface, $y{\in}\partial\Omega.$ The second equality in \eqref{4.1} is immediate from Gauss-Ostrogradsky divergence formula.

 The total density flux through the surface $\partial\Omega$ per unit time has the form
\begin{align}& 
  \label{4.2}  B(f,\Omega):=\int_{\Omega}\int_{\BR^n}v{\cdot}{\bigtriangledown_{x}}f(t,x,v)d{v}dx= \int_{\partial\Omega}\int_{\BR^n}\big(v{\cdot}n(y)\big)f(t,y,v)d{v}d\sigma.
 \end{align}

\begin{Dn} By  $D(\Omega)$, we denote the class of 
functions $f(t,x,v)$ satisfying the Fokker-Planck equation \eqref{1}, inequalities \eqref{3}
and the  condition $A(f,\Omega){\geq}0$ for all $t$ $($i.e., the class of the dissipative solutions $f$$)$.
\end{Dn}

\begin{Dn} By  $C(\Omega)$, we denote the class of 
functions $f(t,x,v)$ satisfying the Fokker-Planck equation \eqref{1}, inequalities \eqref{3}
and the  condition $A(f,\Omega){=}0$ for all $t$ $($i.e., the class of the conservative solutions $f$$)$.
\end{Dn}

It is obvious that $C(\Omega){\subset}D(\Omega).$

\begin{Pn}
Let  $f(t,x,v)$ satisfy \eqref{1} and \eqref{3}, and assume that for all $y\in \partial\Omega$ the equality
 \begin{equation} \label{4.3}
 f(t,y,v)=f(t,y,-v)
 \end{equation}
holds.  Then, $f(t,x,v)\in C(\Omega)$.
 \end{Pn}

\begin{proof}
Taking into account  \eqref{4.3}, we derive
 \begin{equation} \label{4.4}
 \int_{\BR^n}(|v|^{2}/2)\big(v{\cdot}n(y)\big)f(t,x,v){v}d{v}=0.
 \end{equation}
It is immediate from \eqref{4.1} and \eqref{4.4} that  $A(f,\Omega){=}0$.
\end{proof}

\begin{Cy}
Maxwellians  $M_{k}$ given by
\eqref{16} are conservative, that is, $M_k\in C(\Omega)$.
\end{Cy}

The \textbf{bounce-back condition} \eqref{4.3} means that particles arriving with a certain velocity 
to the boundary $\partial\Omega$ will bounce back with an opposite velocity (see \cite[p.16]{Vi1}).

\subsection{Boundedness}
Introducing the function
\begin{equation}
\label{5.1}
s(r)=r\log{r}-\frac{1}{k}(1+{k}r) \log(1+{k}r),\end{equation}
we see that
\begin{equation}\label{5.2} s^{\prime}(r)=\log\big(r/(1+{k}r)\big),\quad
s^{\prime\prime}(r)=\big(r(1+{k}r)\big)^{-1},\end{equation}
and obtain the proposition below.
\begin{Pn}\label{Pn5.1}
The function $s(r)$ is convex on the semi-axis $r \geq 0$ for the case that $k>0$ and on the interval $0{\leq}r<1/|k|$ for the case that $k<0$.
\end{Pn}
We consider functions $g$ such that
\begin{equation}\label{5.5}g\geq 0, \quad 1+kg>0, \quad \frac{1}{2}\int_{\BR^n}|v|^{2}g(v)dv<\infty.
\end{equation}
For $k<1$ and $Z$ given by
\begin{equation}\label{5.3}Z=1/(\E^{|v|^{2}/2}-k),
\end{equation}
the convexity of $s$ implies that
\begin{equation}\label{5.4}s(g)-s(Z){\geq}s^{\prime}(Z)(g-Z).\end{equation}
Using \eqref{5.2}, \eqref{5.3}, \eqref{5.4}  and the equality
\begin{equation}\label{5.6}Z/(1+{k}Z)=\E^{-|v|^{2}/2},\end{equation}
we easily derive
\begin{equation}\label{5.7}-s(g){\leq}\frac{|v|^{2}}{2}g-s(Z)-\frac{|v|^{2}}{2}Z.
\end{equation}
Taking into account the equality
\begin{equation}\label{5.8}1+{k}Z=1/(1-{k}\E^{-|v|^{2}/2}), \end{equation}
we rewrite \eqref{5.7} in the form
\begin{equation}\label{5.9}-s(g){\leq}\frac{|v|^{2}}{2}g
-\frac{1}{k}\log(1-{k}e^{-|v|^{2}/2}).
\end{equation}
The  proposition below follows from \eqref{5.9} and the inequality $\log(1+a){\leq}a$ for $a{\geq}0$.
\begin{Pn}\label{Pr5.2} Assume that $k<0$ and $g$ satisfies  \eqref{5.5}.
Then we have
\begin{equation}
\label{5.11}
-{|v|^{2}}g/2-s(g){\leq}
\E^{-|v|^{2}/2}.
\end{equation}
\end{Pn}
\begin{Rk}\label{Rk5.1.} For the case 
$k=-1$, the
inequality \eqref{5.11} was  derived  in \cite{CPR}.
\end{Rk}
Recall definitions  \eqref{4}, \eqref{9} and \eqref{12} of $S$, $E$ and $F$, respectively,
and note that $S$ is expressed via $\Phi$.
Compare the expression \eqref{5} for $\Phi$ 
 with the expression for $s$.
Thus, we see that  formula \eqref{5.11} yields the following corollary.

\begin{Cy}\label{Cy5.1} Let the conditions of Proposition \ref{Pr5.2} be fulfilled.
Then there exists a positive constant $C$ such that
\begin{equation}\label{5.12}F(g){\leq}C.\end{equation}
\end{Cy}
In a similar way, from \eqref{5.9} and series representation of $\log(1+a)$ we obtain the next corollary.
\begin{Cy}\label{Cy5.2} Let $0<k<1$ and assume that $g$ satisfies  \eqref{5.5}.
Then there exists a positive constant $C$ such that
\begin{equation}\label{5.13} F(g){\leq}C.\end{equation}
\end{Cy}
If either conditions of Corollary \ref{Cy5.1} or Corollary \ref{Cy5.2} hold, we put
\begin{equation}\label{5.13'} \wt C=\inf F(g).\end{equation}
Passing to the limit and using Corollaries \ref{CyB}, \ref{Cy 2.5} and \ref{CyF}, we could show
(under rather general conditions) the existence of the global Maxwellian
$M_{\wt C}$ such that
\begin{align}&
\label{7.9} F(M_{\wt C})=\wt C.
\end{align}
Later we assume that $M_{\wt C}$ exists.
\section{Lyapunov functional}
\label{Lya}
The Lyapunov functional for equation \eqref{1} has the form
 \begin{equation}\label{7.1} \wt G(f)=\wt C-F(f),\end{equation}
 where
 $F(f)$ is defined by \eqref{12}. Recall that, for the generalization $G(f)$ (given by \eqref{18}) of the Kullback-Leibler distance,
 we always assumed that $\rho(f)$ is fixed at the moment $t$ (i.e., $\rho(f,t) \equiv \wt \rho$).
 We do not assume this in the present section. In other words, we consider $\wt G$ on a wider set
 of solutions.
 Clearly, under conditions of Corollaries \ref{Cy5.1} or \ref{Cy5.2}, the inequality
\begin{equation}
\label{7.2} \wt G(f){\geq}0
\end{equation}
 holds. Thus, if a stationary solution $M_{\wt C}$ of \eqref{1} satisfies the equality \eqref{7.9}, 
the value $\wt G(f)=F(M_{\wt C})-F(f)$ at $t$ may be considered as a distance  between $M_{\wt C}$ and $f$ at the time $t$.

In this section we substitute condition \eqref{3} by a stronger condition with the strict inequalities
\begin{equation}\label{3'} 
f(t,x,v) > 0, \quad 1+kf(t,x,v) > 0 \quad (k \in \BR).
\end{equation}
 Hence, we can consider equation \eqref{1} in the form 
  \begin{equation}\label{7.3}
 \frac{\partial{f}}{\partial{t}}=-v \cdot \bigtriangledown_{x}{f}+
 \dv_{v}\left(f(1+kf){\bigtriangledown}_{v}\left(\log\frac{f}{1+kf}+1+|v|^{2}/2\right)\right).\end{equation}
 We multiply \eqref{7.3} { by $\left(\log\frac{f}{1+kf}+1+|v|^{2}/2\right)$ and integrate over $\BR^n$ (applying also integration by parts
 with respect to the  variables $v_i$) 
 and $\Omega.$
Under natural assumptions on the decay of $f$ and $\p f/\p v_i$ at infinity (so that the values at $\pm \infty$ disappear in the formulas
for the integration by parts),  we obtain
\begin{align}& \int_{\Omega}\int_{\BR^n}\left(\log\frac{f}{1+kf}+1+|v|^{2}/2\right)\frac{\partial{f}}{\partial{t}}dvdx {\nonumber} \\
 &=- \int_{\Omega}\int_{\BR^n}\left(\log\frac{f}{1+kf}+1+|v|^{2}/2\right)v \cdot \bigtriangledown_{x}{f}dvdx \label{7.4}\\ & \quad
 {\nonumber} -\int_{\Omega}\int_{\BR^n}f(1+kf)\left|{\bigtriangledown}_{v}\left(\log\frac{f}{1+kf}+1+v^{2}/2\right)\right|^{2}dvdx.
 \end{align} 
 Let us introduce the total flow of the entropy across the boundary $\Omega$:
 \begin{equation}\label{7.5}U(f, \Omega)=-\int_{\Omega}\int_{\BR^n}(\log\frac{f}{1+kf}+1)v \cdot \bigtriangledown_{x}{f}dvdx=
-\int_{\Omega}\int_{\BR^n}v\cdot \bigtriangledown_{x}\Phi(f)dvdx,
\end{equation}
where $\Phi$ is defined by \eqref{5'} and \eqref{5}.
 Using Gauss-Ostrogradsky formula, we also have 
\begin{equation}\label{7.5'}U(f, \Omega)=-
\int_{\partial\Omega} \int_{\BR^n}\big(v{\cdot}n(y)\big)\Phi(f)dvd\sigma .\end{equation}
It follows from \eqref{7.4} and \eqref{7.5}  that
\begin{equation}\label{7.6}\frac{d\wt G}{dt}{\leq}U(f,\Omega)-A(f,\Omega),\end{equation}
where the function $A(f,\Omega)$ is defined by relation \eqref{4.1}.
Indeed, relations \eqref{4.1} and \eqref{7.5}  yield the equality
\begin{align}& U(f,\Omega)-A(f,\Omega)=- \int_{\Omega}\int_{\BR^n}\left(\log\frac{f}{1+kf}+1+|v|^{2}/2\right)v \cdot \bigtriangledown_{x}{f}dvdx,
\label{7.6'}
\end{align}
whereas formulas \eqref{3}, \eqref{7.4} and \eqref{7.6'} imply that
\begin{align}&  \label{7.6''} U(f,\Omega)-A(f,\Omega)
\geq \int_{\Omega}\int_{\BR^n}\left(\log\frac{f}{1+kf}+1+|v|^{2}/2\right)\frac{\partial{f}}{\partial{t}}dvdx.
\end{align}
Furthermore, according to \eqref{4}, \eqref{9}, \eqref{12} and \eqref{7.1}, we have
\begin{equation}\nn
\frac{\p\wt G}{\p t}=\frac{\p E}{\p t}-\frac{\p S}{\p t}=\int_{\Omega}\int_{\BR^n}\left(\log\frac{f}{1+kf}+1+|v|^{2}/2\right)\frac{\partial{f}}{\partial{t}}dvdx,
\end{equation}
and so \eqref{7.6} is immediate from \eqref{7.6''}.

Using inequality \eqref{7.6}, we derive the following assertion.
\begin{Tm}\label{Tm7.1}
Assume that $f{\in}C_{0}^{1}$ is a dissipative solution of  \eqref{1}, which satisfies
\eqref{3},  and  that the inequality
\begin{equation}\label{7.8} U(f,\Omega){\leq}0 \end{equation}
holds.  Then the inequality
 $( \p{\wt G}/{\p t}){\leq}0$ is valid.
\end{Tm}
\begin{Cy}
Assume that the  conditions of Theorem \ref{Tm7.1} are fulfilled and
$\wt G(f,t_0)<\delta$. Then, the inequality $\wt G(f,t)<\delta$
holds for all $ t>t_{0}$.
\end{Cy}
Thus, the distance $\wt G$ between a Maxwellian $M_{\wt C}$ satisfying \eqref{7.9} and $f$
satisfying \eqref{7.8} decreases. Taking into account the definition of the Lyapunov stability,
we proved the following important result.
\begin{Tm}\label{Tm7.2} The stationary solution $M_{\wt C}$ is locally  stable $($i.e., Lyapunov stable$)$
in the class of the dissipative solutions $f$ satisfying inequalities \eqref{3} and
\eqref{7.8}. \end{Tm}
We note that $M_{\wt C}$ does not depend on $x$, and so definitions \eqref{4.1} and \eqref{7.5} imply
that $A(M_{\wt C}, \, \Omega)=U(M_{\wt C}, \, \Omega)=0$.
 
\begin{Rk}The earlier results on the local stability for the Fokker-Planck equation
$($see \cite{CPR}$)$ were obtained  only for   the spatially homogeneous case. Theorems \ref{Tm7.1} and \ref{Tm7.2} have no analogs also 
in the case of
the Boltzmann equation.\end{Rk}

{\bf Concluding remark.}
Following \cite{LAS1, LAS2, LAS3, LASPLet, LAS4}, we compare classical and quantum results, that is,
determinate and probabilistic cases.
As in the case of the Boltzmann equation \cite{LASPLet, LAS4},
we use a special extremal principle, which is based on the ideas of  game theory. (We note that extremal principles  remain central in modern physics.) 
The "players" in the game described by the Fokker-Planck equation
are the total energy $E$ and entropy $S$, the "gain" in the game is the functional $F$ and the strategy in the game  is determinate in the classical case
and probabilistic in the quantum case. It is of interest that the inequalities $E_{q,B}<E_{c}<E_{q,F}$, $S_{q,B}<S_{c}<S_{q,F}$ (for $n>2$), and $F_{q,F}<F_{c}<F_{q,B}$
(see \eqref{3.34}) hold for small $k$ in our game. 

\noindent{\bf Acknowledgments.}
 {The research of A.L. Sakhnovich was supported by the
Austrian Science Fund (FWF) under Grant  No. P24301.}

\end{document}